\documentstyle[psfig]{JKTR-amslatex2}
\def\wx{w_{x}}
\def\wy{w_{y}}

\def\S2{S^{2}}

\newcommand{\illu}[3]{
  \begin{center}
  \mbox{
    \rule{0mm}{#2mm}
    \psfig{file=#3}
  }
  \end{center}
}
\newcommand{\smillu}[3]{
  \mbox{
    \rule{0mm}{#2mm}
    \psfig{file=#3}
}}

\begin{document}

\title{NULLIFICATION WRITHE AND\\CHIRALITY OF ALTERNATING LINKS}

\catcode`\@=11

\author{CORINNE CERF\thanks{The author is 
``Charg\'e de Recherches" of the Belgian ``Fonds National de la
Recherche Scientifique". E-mail: {\tt ccerf@ulb.ac.be}}\\
{\it D\'epartement de Math\'ematiques, CP 216}\\
{\it Universit\'e Libre de Bruxelles}\\
{\it Boulevard du Triomphe}\\
{\it B-1050 Bruxelles, Belgium}\\
{\it (mailing address)}\\
{\it and}\\
{\it Department of Chemistry}\\
{\it Princeton University}\\
{\it Princeton, NJ 08544, U.S.A.}}

\catcode`\@=\active

\maketitle

\vspace{0.1cm}
Published in {\it J. Knot Theory Ramifications} {\bf 6} (1997) 621-632.

\begin{abstract}
In this paper, we show how to split the writhe of reduced projections of oriented 
alternating links into two parts, called the nullification writhe $\wx$, and the remaining 
writhe $\wy$, such that the sum of these quantities equals the writhe $w$ and each 
quantity remains an invariant of isotopy. The chirality of oriented alternating links can be 
detected by a non-zero $\wx$ or $\wy$, which constitutes an improvement compared to 
the detection of chirality by a non-zero $w$. An interesting corollary is that all oriented 
alternating non-split links with an even number of components are
chiral, a result that also follows  from properties of the Conway
polynomial.
\vskip 1em \noindent{\it Keywords}: Link, chirality,
writhe, nullification writhe, remaining writhe. \end{abstract}

\section{Introduction}

Determining whether knots and links are equivalent or not to their mirror images, i.e., 
whether they are $achiral$ or $chiral$, has been a longstanding question in knot theory.
Tait  pioneered this field in the nineteenth century by developing empirical methods to
this end  \cite{Tai}. Subsequently, more and more sophisticated methods capable of
detecting the  chirality of links (essentially numerical and polynomial invariants) were
elaborated by  numerous mathematicians. A critical review of these methods can be found
in  \cite{Fla91,LiMi95}. 

In 1991, Sola introduced a numerical invariant for oriented alternating links, called 
the {\it nullification number} \cite{Sol91}. This invariant is, however, insensitive to
chirality. We  show here that by a modification of the nullification number, a new
numerical invariant  capable of detecting chirality can be derived. It is called
the {\it nullification writhe}.

\section{Concept of nullification}

Let us begin with some definitions.

\begin{figure}[htp]
\illu{41}{12}{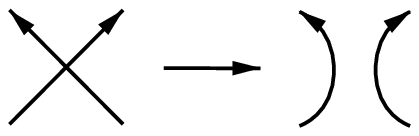}
\caption{}
\label{fig:1}
\end{figure}

\begin{defin}\label{def:null_cross}
The {\it nullification of a crossing} of an oriented link
projected on $\S2$ is the process  described in Fig.~\ref{fig:1},
where\raisebox{-1mm}{\smillu{4}{4}{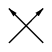}}
means\raisebox{-1mm}{\smillu{4}{4}{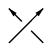}}
or\raisebox{-1mm}{\smillu{4}{4}{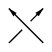}}.  
\end{defin}

\begin{defin}\label{def:null_proj}
The {\it nullification of an oriented link projection} on $S^2$
consists in nullifying a series of
crossings until an unknot (or unlink in the case of a
split link) is reached, while preventing at each step the disconnection
of the link (or link components in the case of a split link).
\end{defin}

Let us observe that when the stage of the unknot (or unlink) is
reached, all remaining crossings are nugatory crossings, so that
nullifying any one of them would disconnect the unknot (or add a
disconnected component to the unlink). Conversely, any projection
composed of nugatory crossings only, is necessarily an unknot or an
unlink.

We will call the set of crossings so nullified the {\it nullification 
set}. Sola proves that the cardinality of this 
set, called the {\it nullification number} and denoted by $o$, is
well defined and invariant for reduced alternating link
projections \cite{Sol91}:

\begin{equation}
o(K) = n(K)-s(K)+1
\end{equation}

\noindent where $n(K)$ is the number of crossings of a reduced alternating
projection of link $K$ and $s(K)$ is the number of Seifert circles of 
a reduced alternating projection of link $K$ (obtained by
simultaneously nullifying all crossings).

We would like to prove that the writhe of the nullification set (the {\it
nullification writhe}) is also well defined and invariant for reduced
alternating link projections. 

\begin{defin}
The {\it nullification writhe} $\wx$ of an oriented link projection
on $S^2$ is the sum of the signs of the crossings extracted during the
nullification of the projection.
\end{defin}

\begin{defin}
The {\it remaining writhe} $\wy$ of an oriented link projection on
$S^2$ is the sum of the signs of the crossings of the unknot (or
unlink in the case of a split link) remaining after nullification of
the projection.
\end{defin}

It follows immediately that for a given projection of a link $K$,
\begin{equation}\label{w}
w(K) = \wx(K) + \wy(K).
\end{equation}

\begin{figure}[htp]
\illu{70}{34}{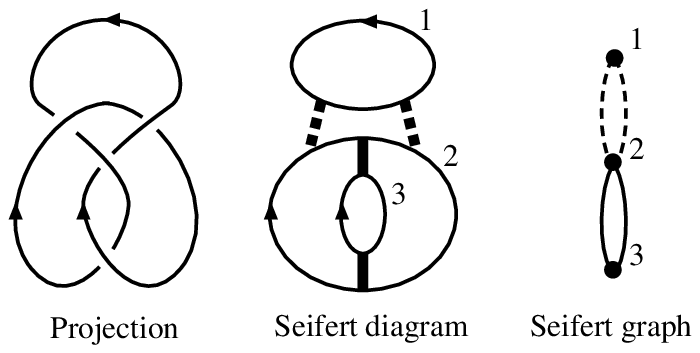}
\caption{}
\label{fig:proj_diag_graph}
\end{figure}

We will make use of two other descriptions of the link projection.
The first one is called {\it Seifert diagram}. It is obtained by
simultaneously nullifying all crossings (this gives rise to a
collection of Seifert circles) and by replacing each crossing of the link
projection with an arm, called $connection$, connecting both involved
Seifert circles. Positive and negative connections are distinguished
and represented by a solid bar and a  dashed bar, respectively. The
second description is called {\it Seifert graph}. It is obtained from
the Seifert diagram by crushing each Seifert circle to a point,
that becomes a vertex of the graph. Connections become positive and
negative edges of the graph. An example is shown in 
Fig.~\ref{fig:proj_diag_graph}.

The process of nullification of a link projection transposed in the
Seifert graph representation consists in removing the maximum number
of edges without disconnecting the graph (or without adding any
additional disconnected component to the graph, in the case of a
split link). This process leaves a spanning tree (respectively, a
spanning forest). The starting graph possesses $s$ vertices (the
number of Seifert circles) and $n$ edges
(the number of crossings). It is well known that the number of
edges of a spanning forest for a split graph with $s$ vertices and $k$
components is $s-k$ (see \cite{Wil85} p. 45 for example). It follows
that the nullification number, i.e., the number of edges that have
been removed, is $n-s+k$, which confirms Sola's result (Sola assumes
implicitely that the link is non-split, so $k=1$).

\section{Well-definedness}

In order to prove that $\wx$ and $\wy$ are well defined, we will need 
the following lemma:

\begin{lem}\label{lem:sides}
In the Seifert diagram corresponding to an alternating link
projection, all connections on the same side of a Seifert circle
have the same sign, and all connections on opposite sides of a
Seifert circle have opposite signs.
\end{lem}

\begin{figure}[htp]
\illu{60}{28}{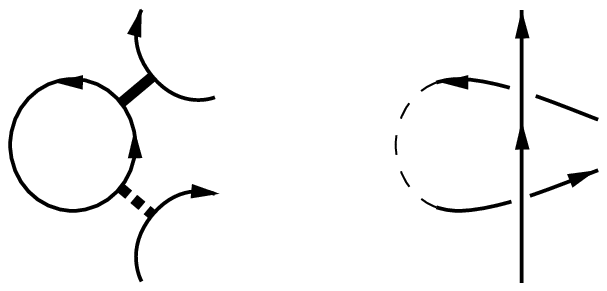}
\caption{}
\label{fig:same_side}
\end{figure}

\begin{figure}[htp]
\illu{62}{27}{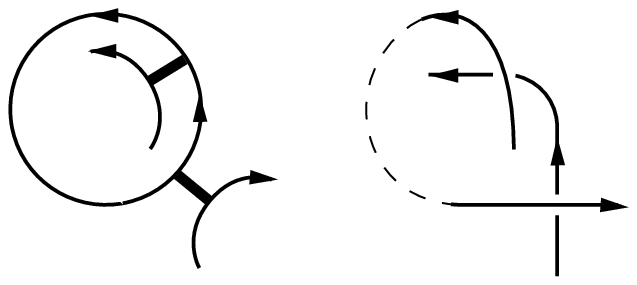}
\caption{}
\label{fig:opposite_sides}
\end{figure}

\begin{pf}
Two adjacent connections on the same side of a Seifert circle can
never have opposite signs, because the corresponding link projection
would be non-alternating (Fig.~\ref{fig:same_side}). Two
adjacent connections, one on one side and one on the other side of a
Seifert circle can never have the same sign, because the corresponding
link projection would be non-alternating
(Fig.~\ref{fig:opposite_sides}).
\end{pf}

We can now prove

\begin{prop}\label{prop:well_defined}
The nullification writhe and the remaining writhe of a reduced
projection of an oriented alternating link are independent of
the choice of the nullification set. 
\end{prop}

\begin{pf}
We know that any two nullification sets of a reduced projection of an
oriented alternating link have the same cardinality. Do they have the
same writhe (nullification writhe)? In the Seifert graph description,
the nullification set contains all edges to remove in order to
obtain a spanning tree (or spanning forest if the link is split), that is, 
all edges but one between each pair of vertices linked by multiple edges,
and one edge from each cycle of edges. By Lemma~\ref{lem:sides}, all
multiple edges in a Seifert graph (corresponding to multiple
connections in a Seifert diagram) have the same sign, so keeping any one of 
them will add the same contribution to the nullification
writhe. As far as cycles are concerned, because of Lemma~\ref{lem:sides} 
and of the Jordan Curve
Theorem (see \cite{Arm83} for example), all cycles consist of edges
having the same sign, so removing any edge will add the same
contribution to the nullification writhe. The nullification writhe 
is thus independent of the choice of
the nullification set. 
From Eq.~\ref{w} (or simply for complementarity reasons), it follows that
the remaining writhe is also independent of the choice of the
nullification set. 
\end{pf}

\section{Invariance}

Let us begin by proving the following lemma that will simplify the proof 
of invariance:

\begin{lem}\label{lem:any_crossing}
Given any crossing of a reduced projection of an oriented alternating
link, there exists a nullification set containing this crossing.
\end{lem}

\begin{pf}
The fact that the projection is reduced means that there is no nugatory crossing in
it. In the Seifert diagram description, this means that there is no
connection whose removal would disconnect the diagram (or
add a disconnected component in the case of a split link). In the
nullification process, we may thus choose to begin by the extraction
of any crossing, that will become part of the nullification set.
\end{pf}

\begin{prop}\label{prop:invariant}
Any two reduced projections of an oriented alternating link have
the same nullification writhe and the same remaining writhe.
\end{prop}

\begin{pf}

As asserted by the Tait Flyping Conjecture \cite{Tai} proved by
Menasco and  Thistlethwaite \cite{MeTh91}, any two reduced projections
of an alternating link are  related by a sequence of moves called
$flypes$. It remains to prove that two oriented  reduced alternating
projections related by a flype have the same nullification writhe and
the same remaining writhe.

\begin{figure}[p]
\illu{85}{19}{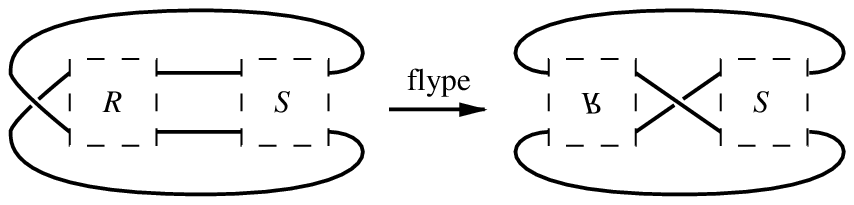}
\caption{}
\label{fig:final05}
\end{figure}

Two projections related by a flype can be represented as in Fig.~\ref{fig:final05}, where $R$ and 
$S$ are tangles, i.e., parts of a link projection with four emerging arcs. The depicted 
crossing may be positive or negative. The sign of the crossing is unchanged by the 
flyping operation. Let us now represent this operation in the description of Seifert 
diagrams. Depending on the orientation, two situations occur, shown in
Figs.~\ref{fig:final06}~and~\ref{fig:final08}. The highlighted connection is
represented as a positive connection, but it could also be a
negative connection, for the same reason as above.

\begin{figure}[p]
\illu{85}{19}{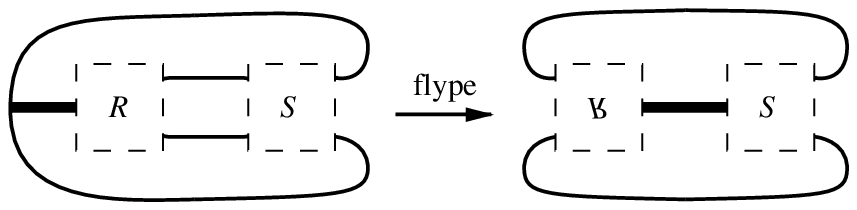}
\caption{}
\label{fig:final06}
\end{figure}

\begin{figure}[p]
\illu{85}{19}{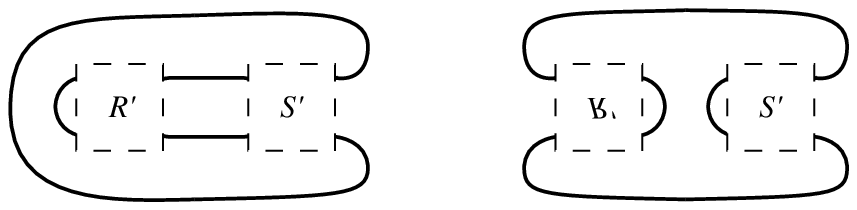}
\caption{}
\label{fig:connected_sum2}
\end{figure}

\begin{figure}[p]
\illu{86}{19}{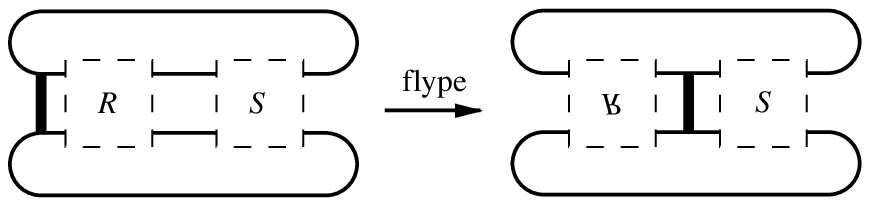}
\caption{}
\label{fig:final08}
\end{figure}

\begin{figure}[p]
\illu{86}{19}{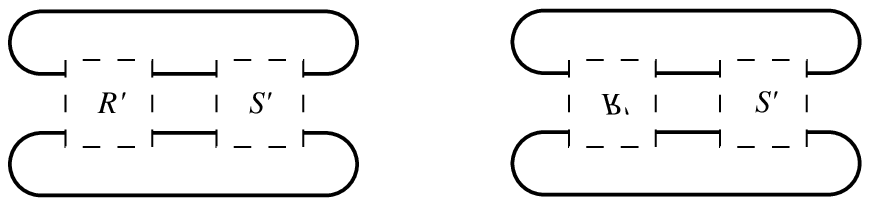}
\caption{}
\label{fig:connected_sum3}
\end{figure}

Let us analyze the first situation (Fig.~\ref{fig:final06}). The flype transforms
each connection into a connection of the same sign. 
We now nullify both diagrams of Fig.~\ref{fig:final06}. For the left-hand diagram, 
we choose a nullification set containing the highlighted
connection (which is legitimate because of
Lemma~\ref{lem:any_crossing}). It remains a connected sum of two tangles 
$R'$ and $S'$ (Fig.~\ref{fig:connected_sum2}, left). Knowing that it is an unknot 
(or an unlink), tangles $R'$ and $S'$ are necessarily
unknots (or unlinks). Then, for the diagram on the right, 
we choose as nullification set the image of the left-hand nullification 
set. Is this really a 
nullification set? What is left is a
connected sum of the same tangles $R'$ and $S'$ with a different relative
orientation (Fig.~\ref{fig:connected_sum2}, right). Since $R'$ and $S'$ are
unknots (or unlinks), the right-hand diagram is an unknot (or unlink). 
This completes the proof for the first situation.

Let us now analyze the second situation (Fig.~\ref{fig:final08}). Again, each 
connection is transformed into a connection of the same sign.
We nullify both diagrams (Fig.~\ref{fig:connected_sum3}). For the left-hand diagram, we
choose a nullification set containing the highlighted
connection (which is legitimate because of
Lemma~\ref{lem:any_crossing}). For the right-hand diagram, we use the
image of this set.

\begin{figure}[p]
\illu{86}{19}{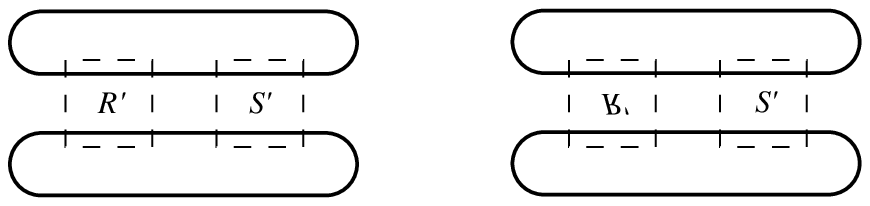}
\caption{}
\label{fig:caseI}
\end{figure}

\begin{figure}[p]
\illu{86}{19}{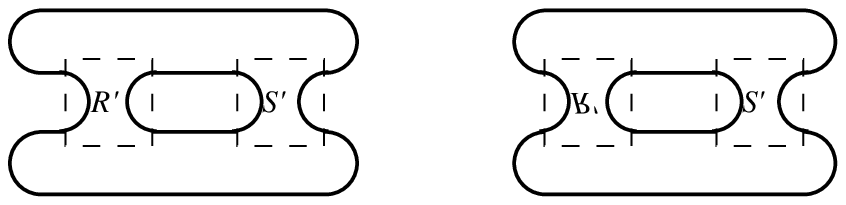}
\caption{}
\label{fig:caseII}
\end{figure}

\begin{figure}[p]
\illu{86}{19}{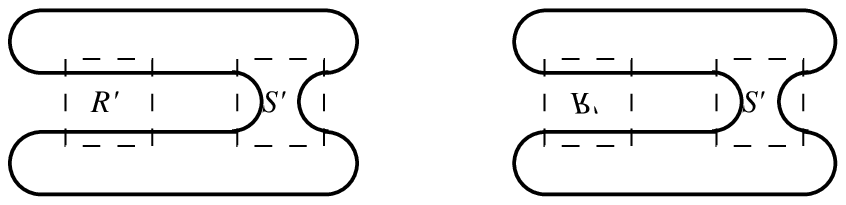}
\caption{}
\label{fig:caseIII}
\end{figure}

\begin{figure}[p]
\illu{86}{19}{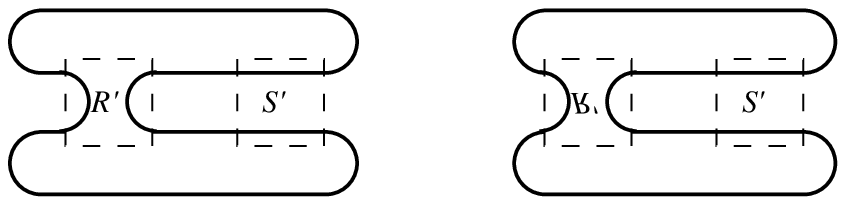}
\caption{}
\label{fig:caseIV}
\end{figure}

Four subcases arise, represented in Figs.~\ref{fig:caseI}-\ref{fig:caseIV}.
In the first two subcases
(Figs.~\ref{fig:caseI}-\ref{fig:caseII}), we know that at the end of the
nullification process, there must remain only one connection (direct or
indirect) between both displayed circles. One of \{$R',S'$\} may
thus be split into 2 tangles (horizontal splitting in the first subcase,
vertical splitting in the second subcase). The diagrams become connected
sums of 3 tangles. Because the left-hand diagrams are unknots (or
unlinks), the 3 tangles are necessarily unknots (or
unlinks). Therefore, the right-hand diagrams, that are connected sums
of the same 3 tangles with different orientations, also are unknots (or 
unlinks). The image of a nullification set is thus a nullification set.

Let us finally consider the last two subcases 
(Figs.~\ref{fig:caseIII}-\ref{fig:caseIV}). By the definition of a 
connection, there cannot exist a connection connecting a circle to itself. 
Therefore, both $R'$ and $S'$ may be split into 2 tangles (in each subcase,
one horizontal and one vertical splitting). Let us split one tangle only.
This leads again to connected sums of 3 tangles. Because the left-hand 
diagrams are unknots (or
unlinks), the 3 tangles are necessarily unknots (or
unlinks). Therefore, the right-hand diagrams, that are connected sums
of the same 3 tangles with different orientations, also are unknots (or 
unlinks). The image af a nullification set is thus a nullification set.

This  completes the proof of Proposition \ref{prop:invariant}.   
\end{pf}

\section{Detection of Chirality}

As was known for long, an easy way to detect the chirality of an
oriented alternating link
is to look at the  writhe of a reduced projection of the link. If it is 
non-zero, then the link is chiral.
Here is a very simple proof of this fact.

\begin{prop}
If an oriented alternating link $K$ is achiral, then the writhe of a
reduced projection of $K$ is equal to zero.
\end{prop}

\begin{pf}
Flypes preserve the writhe of reduced projections of oriented
alternating links. If a link is achiral, the proven 
Tait Flyping Conjecture says that its mirror image can be obtained by a
series of flypes. Thus a projection and its mirror image have
at the same time identical writhes (because of the Tait Flyping
Conjecture) and opposite writhes (because in the mirror image, all
crossings have changed sign). The writhe of the link must then be
equal to zero.
\end{pf}

The chirality of an alternating link can thus be detected by a
non-zero writhe. Unfortunately, this is not a universal way to
detect chirality:  there exist many chiral links whose 
writhe is equal to zero too.

The same problem of zero-writhe chiral links has been encountered in
a previous study  aimed at partitioning chiral knots and links into
$D$ and $L$ classes, and was overcome  by defining a property called
{\it writhe profile} \cite{LiMi94,LiCeMi96}. This property is of no  help
here, because it is a ``chirality-classifier" and not a
``chirality-detector". The chirality  of the link is a prerequisite in
order to apply the method. In contrast, the new invariants  introduced
in this paper, namely the nullification writhe $\wx$ and the remaining
writhe $\wy$, are capable of detecting chirality, as will be proved below. 

\begin{prop}\label{prop:mirror}
If $K$ is an oriented alternating link represented by a reduced
projection and if $K^*$ is its mirror image, then $\wx(K)=-\wx(K^*)$ and 
$\wy(K)=-\wy(K^*)$.
\end{prop}

\begin{pf}
Given a nullification set for $K$, let us choose for $K^*$ the
corresponding nullification set, i.e., the crossings located at the
same position and having the opposite sign. The writhes of these two
sets are obviously opposite, so
$\wx(K)=-\wx(K^*)$. Similarly, the sums of the signs of the remaining
crossings in $K$ and $K^*$ are opposite, so
$\wy(K)=-\wy(K^*)$.
\end{pf}

\begin{cor}
If an oriented alternating link $K$ is achiral, then $\wx(K)$ and
$\wy(K)$, computed on a reduced projection of $K$, are equal to zero.
\end{cor}

\begin{pf}
As we proved in Proposition~\ref{prop:invariant}, flypes preserve 
$\wx$ and $\wy$ for reduced projections of oriented
alternating links. If a link $K$ is achiral, the Tait Flyping
Conjecture says that its mirror image, $K^*$, can be obtained by a
series of flypes. Thus $K$ and $K^*$ have
at the same time identical $\wx$ (respectively, $\wy$) because of the
Tait Flyping Conjecture, and opposite $\wx$ (respectively, $\wy$)
because of Proposition~\ref{prop:mirror}. It
implies that $\wx(K)=0$ and $\wy(K)=0$. 
\end{pf}

We have thus the following implication: if $\wx(K) \neq 0$ or $\wy(K)
\neq  0$, then $K$ is chiral. What are the consequences of
this? For each oriented  alternating link represented by a reduced
projection, the writhe $w$ is split into two parts $\wx$ and $\wy$.
If $w$ is different from zero, $\wx$ and/or $\wy$ are different from 
zero and the link is chiral. In this case, $\wx$  and $\wy$ are thus
no better than $w$ in  detecting chirality. In contrast, if
$w$ equals zero but the link is chiral, then in some  cases $\wx$  and
$\wy$ might be different from zero ($\wy = -\wx$), which would 
give a way to detect the chirality of the link even if it is not
detectable by $w$.

We have applied this procedure to all chiral oriented alternating ``classical" knots and
links,  i.e., prime knots with up to 10 crossings and non-split prime links with up to 9
crossings and  4 components, that have a writhe of zero. The results are shown in
Tables 1 and 2 and  confirm that chirality is detected in several cases where it was not
detectable by the writhe only.

\begin{table}[t]
\renewcommand{\arraystretch}{2}
\caption{Chiral oriented alternating prime knots
with up to 10 crossings, and $w = 0$ ($\wy = -\wx$).}
\small{\begin{tabular}[t]{@{} l c @{ } c @{\hspace{2mm}} c @{\hspace{2mm}} c
@{\hspace{2mm}} c  @{\hspace{2mm}} c @{\hspace{2mm}} 
c @{  } c @{\hspace{2mm}} c @{\hspace{2mm}} c @{\hspace{2mm}} c @{\hspace{2mm}} c 
@{\hspace{2mm}} c @{\hspace{2mm}} c @{\hspace{2mm}} c @{}} \hline
knot\footnotemark[1]$^{,}$\footnotemark[2] & \hfill & $8_{4}$ & $10_{15}$ & $10_{19}$ & $10_{31}$ & $10_{42}$ & 
$10_{48}$ & $10_{52}$ & $10_{54}$ & $10_{71}$ & $10_{91}$ & $10_{93}$ & $10_{104}$ & 
$10_{107}$ & $10_{108}$ \\ \hline
$\wx$ && $-2$ & +2 & $-2$ & 0 & 0 & 0 & +2 & +2 & 0 & 0 & $-2$ & 0 & 0 & +2 \\
\hline
\end{tabular}}
\end{table}

\begin{table}[t]
\renewcommand{\arraystretch}{2}
\caption{Chiral oriented alternating non-split prime links
with up to 9 crossings and 4 components, and $w = 0$ ($\wy = -\wx$).}
\small{\begin{tabular}[t]{@{} l c c @{\hspace{2mm}} c @{\hspace{2mm}} c
@{\hspace{2mm}} c  @{\hspace{2mm}} c
@{\hspace{2mm}} c @{\hspace{2mm}} c @{\hspace{2mm}} c @{\hspace{2mm}} c @{}} \hline
link\footnotemark[3] & \hfill & $8^{2}_{6}\!+\!-$ & 
$8^{2}_{7}\!+\!+$ & $8^{2}_{10}\!+\!+$ & $8^{2}_{11}\!+\!-$ & $8^{2}_{14}\!+\!-$ & 
$8^{3}_{1}\!+\!-\!+$ & $8^{3}_{2}\!+\!-\!-$ & $8^{4}_{1}\!+\!+\!+\!-$ & 
$8^{4}_{1}\!+\!+\!-\!-$ \\ \hline
$\wx$ && $-1$ & +1 & $-1$ & +1 & $-1$ & $-2$ & 0 & $-1$ & $-1$ \\ \hline
\end{tabular}}
\end{table}

\footnotetext[1]{The enantiomorphs are those represented in \cite{Rol76}, oriented either 
way.}
\footnotetext[2]{This list should also contain the non-invertible oriented knots whose 
non-oriented version is achiral. Because of non-invertibility, these knots are chiral when 
oriented. They have been omitted since the nullification writhe is insensitive to 
non-invertibility.}
\footnotetext[3]{The enantiomorphs are those represented in \cite{DoHo91}, using the 
same convention of labeling and orientation of components.}

\section{Links with an even number of components}

This last section will be devoted to the case of oriented alternating
non-split links with an even number of components, for which the
properties of the nullification writhe lead to a very nice proposition:

\begin{prop}\label{prop:even}
All oriented alternating non-split links with an even number of
components are chiral.
\end{prop}

\begin{figure}[htp]
\illu{95}{61}{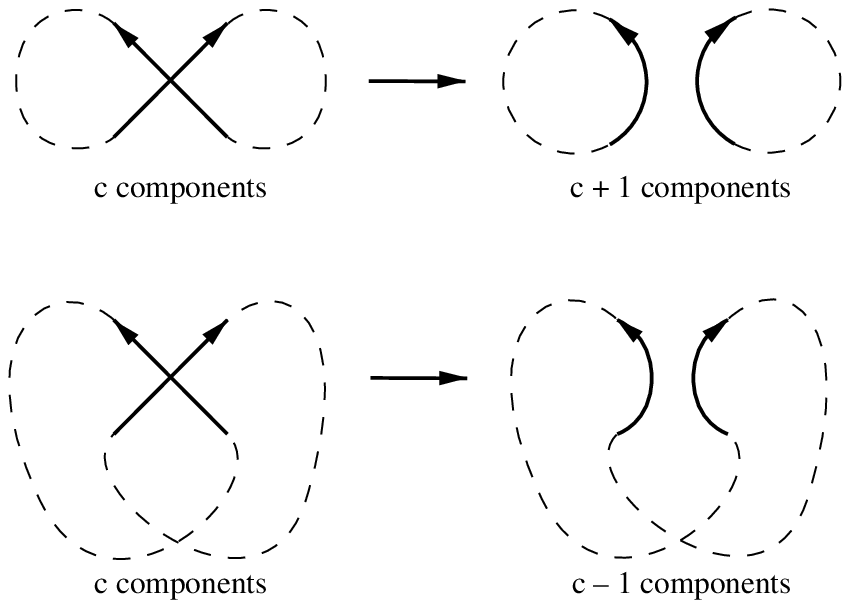}
\caption{}
\label{fig:final14}
\end{figure}

\begin{pf}
Let us nullify an oriented alternating non-split link $K$ with an even
number of components. Each step in the nullification process changes the parity
of the number of components of the link, as can be easily understood
by looking at Fig.~\ref{fig:final14} where\raisebox{-1mm}
{\smillu{4}{4}{ccerf-cross.eps}}
means\raisebox{-1mm}{\smillu{4}{4}{ccerf-crossP.eps}}
or\raisebox{-1mm}{\smillu{4}{4}{ccerf-crossN.eps}}. If the  two arrows
of a crossing are part of the same component, the nullification of the
crossing  increases by one the number of components. If the two arrows
belong to two different  components, the nullification of the crossing
decreases by one the number of components. On the other hand, at the
end of the nullification process, we get an unknot (with one component)
because $K$ is non-split. This implies that the nullification process
contains an odd number of steps. Since each step contributes $\pm1$ 
to $\wx$, $\wx$ can never be equal to zero, which implies that $K$ is
chiral.
\end{pf}

John Conway pointed out that this also follows from properties of his 
$\nabla$ polynomial 
\cite{Con70}, which satisfies

\begin{equation}
\nabla(K^*) = (-1)^{c+1} \nabla(K),
\end{equation}

\noindent where $K^*$ is the mirror image of $K$ and $c$ is the number of components of
$K$.  This implies that an oriented link with an even number of components can only be
achiral if its $\nabla$Êpolynomial vanishes identically. This includes our corollary,
since it is  known that all oriented alternating non-split links
have non-zero Alexander polynomials  \cite{Cro59} and thus non-zero
$\nabla$Êpolynomials ($\nabla$Êpolynomials are  normalized Alexander
polynomials).

\section*{Acknowledgements}
Anna-Barbara Berger, Ines Stassen, and Urs Burri are gratefully
acknowledged for a  crucial discussion held in Bern. This work has
benefited from many fruitful conversations with Chengzhi
Liang and Kurt Mislow. Also Erica Flapan is heartily  thanked for her
critical reading of the manuscript and her invaluable comments.
Finally, I  would like to express my gratitude to John Conway for his
helpful suggestions, and to one of the referees
for his/her detailed and illuminating review. This  work was supported
in part by Grant No CHE-9401774 to K. Mislow.

\end{document}